\documentclass[reqno,11pt]{article}
\usepackage{amsmath,amssymb,euscript,theorem}
\newcommand{\la}{\lambda}
\newcommand{\al}{\alpha}

\newcommand{\fy}{\varphi}
\newcommand{\pd}{\partial}
\newcommand{\p}{\partial}

\newcommand{\R}{\mathbb{R}}
\newcommand{\C}{\mathbb{C}}
\newcommand{\N}{\mathbb{N}}
\newcommand{\Z}{\mathbb{Z}}
\newcommand{\F}{\mathcal{F}}
\renewcommand{\S}{\mathcal{S}}
\renewcommand{\Re}{\mathop{\mathrm{Re}}}
\renewcommand{\Im}{\mathop{\mathrm{Im}}}

\newcommand{\donothing}[1]{}

\numberwithin{equation}{section}
\newtheorem{theorem}{Theorem}[section]

\newtheorem{lemma}[theorem]{Lemma}

\newcommand{\norm}[1]{\left \| #1 \right \|}

\newcommand{\lec}{{\ \lesssim \ }}
\newcommand{\gec}{{\ \gtrsim \ }}
\newcommand{\EQAL}[1]{\begin{equation} \begin{split} #1
 \end{split} \end{equation}}
\newcommand{\EQ}[1]{\begin{equation} #1 \end{equation}}

\newcommand{\de}{\delta}
\newcommand{\si}{\sigma}

\newcommand{\be}{\beta}
\newcommand{\na}{\nabla}
\renewcommand{\th}{\theta}

\newcommand{\supp}{\operatorname{supp}}

\renewcommand{\P}{\mathcal{P}}

\newcommand{\myproof}{{\bf Proof.} \ }
\newcommand{\myendproof}{\quad $\square$}

\begin{document}

\title{Scattering for the Gross-Pitaevskii equation}
\author{Stephen Gustafson, \quad Kenji Nakanishi, \quad Tai-Peng Tsai}
\maketitle
\begin{abstract}
We investigate the asymptotic behavior at time infinity of solutions close to a non-zero constant equilibrium for the Gross-Pitaevskii (or Ginzburg-Landau-Schr\"odinger) equation. We prove that, in dimensions larger than 3, small perturbations can be approximated at time infinity by the linearized evolution, and the wave operators are homeomorphic around 0 in certain Sobolev spaces.
\end{abstract}

\section{Introduction}

The cubic Schr\"odinger equation for 
$\fy(t,x):\R\times\R^d\to\C$, 
\EQ{ \label{NLS}
 i \p_t \fy = -\Delta \fy + |\fy|^2\fy,}
has numerous physical applications.
Noting that the single mode uniform oscillation $e^{-it}$ 
is trivially a solution, and writing $\fy = e^{-it}\psi$, we arrive at 
\begin{equation}\label{GP}
  i \pd_t \psi = - \Delta \psi + (|\psi|^2-1)\psi, \qquad
  \psi(0,\cdot) = \psi_0,
\end{equation}
the {\em Ginzburg-Landau-Schr\"odinger} or {\em Gross-Pitaevskii}
equation, which models the dynamics of Bose-Einstein condensates, 
superfluids, and superconductors (see \cite{FS} for a review).
In these physical settings, $\fy = e^{-it}$ (or $\psi = 1$) 
corresponds to a stationary, constant-density condensate. 
Typically, one is interested in local perturbations of such an equilibrium,
in which case the natural boundary condition as $|x| \to \infty$
should {\em not} be the condition $\psi \to 0$ usually 
imposed when studying~(\ref{NLS}), but rather
\EQ{ \label{BC}
 |\psi| \to 1 \quad \mbox{ as } \quad |x| \to \infty.}

Equation~\eqref{GP} with boundary condition~\eqref{BC} 
possesses, in dimension $d=2$, well-known static,
spatially localized, topologically non-trivial {\em vortex} 
solutions of the form 
$\psi(x) = e^{i n \cdot arg(x)} f(|x|)$, $n \in \Z$.
There are also traveling-wave solutions,  
$\psi(t,x) = v(x_1 - ct, x_2, \cdots, x_d)$,
with $v$ resembling a pair of vortices ($d=2$), 
or a ring of vorticity ($d \geq 3$): \cite{BS,BOS,Ch,Gr}.
For various Ginzburg-Landau Schr\"odinger and wave-type equations 
(such as~\eqref{GP}), the dynamics of vortex points ($d=2$), 
filaments ($d=3$), etc., has been extensively studied 
in certain asymptotic limits:
see, e.g., \cite{CJ,LX,J,L,J2,Sp,GS} and references therein.
(In addition, there is a large literature on vortex dynamics in  
parabolic versions of~\eqref{GP}.)

Very little, however, seems to be known about the long-time 
behaviour of solutions of~\eqref{GP}--~\eqref{BC}, and, 
in particular, about the 
stability of the various interesting vortex-type solutions
mentioned above (however see, eg., \cite{JPR,OS} for numerical and 
formal asymptotic studies of vortex traveling wave stability, and 
\cite{Gu} for vortex stability/instability in the gauge-invariant 
wave-equation version of~\eqref{GP}).
 
As a first step in this direction, we consider
here a much simpler problem: asymptotic stability
of the {\em vacuum} solution
\begin{equation}
\label{vacuum}
  \psi(t,x) \equiv e^{i\theta_0} = \text{constant} 
  \in S^1.
\end{equation}
By the gauge change $\psi\to e^{-i\th_0}\psi$, we may assume $\th_0=0$. 
Further writing $\psi=1+u$, we have 
\EQ{ \label{GLS2}
\begin{split}
&i \pd_t u + \Delta u - 2\Re u = F(u), \\
&F(u) := (u + 2\bar u+ |u|^2)u.
\end{split}} 
The main goal of this paper is to investigate the asymptotic stability of the vacuum in terms of the scattering theory around the constant solution~\eqref{vacuum} of \eqref{GP}, or the $0$ solution of \eqref{GLS2}. 

The first step is to introduce the following change of variable in order to represent the free (linear) evolution as a unitary group:
\EQ{
 u \mapsto v := V^{-1}u := U^{-1}\Re u + i \Im u, \quad U:=\sqrt{-\Delta(2-\Delta)^{-1}}.}
Then $v$ satisfies the equation
\EQ{ \label{eq v}
 i \p_t v - \sqrt{-\Delta(2-\Delta)} v = -i V^{-1} i F(Vv).}
The singularity of $U^{-1}$ at the Fourier origin seems a natural effect of the interaction with the nonzero constant background in \eqref{GP}. Due to this singularity, together with a quadratic nonlinearity which is not gauge invariant, the scattering problem for \eqref{GLS2} turns out to be much more difficult than the original NLS \eqref{NLS} with spatially decaying data (i.e. $\fy\in L^2$). Thus we have obtained so far the scattering only for higher dimensions $d\ge 4$. Global wellposedness of \eqref{GLS2} in the energy space $u\in H^1(\R^d)$ was proved in \cite{BS} for $d=2,3$. It is based on a priori bounds of the energy and so does not work for $d\ge 4$, where the cubic nonlinearity is either energy-critical or super critical. Our theorem gives global wellposedness for $d\ge 4$ and small data as a by-product. 

\begin{theorem}[Stability and scattering]\label{th:1-2}
Suppose that $d \ge 4$ and $|\si| \le \frac{d-3}{2}-\frac 1d$. 
If $U^\si V^{-1}u(0)$ is sufficiently small in $H^{d/2-1}(\R^d)$, then $U^\si V^{-1}u(t)$ remains small in $H^{d/2-1}(\R^d)$ for all $t\in\R$. Moreover, there exist $v_\pm\in U^{-\si}H^{d/2-1}(\R^d)$ such that 
\EQ{
\norm{U^{\si}(V^{-1}u(t) - e^{-iHt}v_\pm)}_{H^{d/2-1}(\R^d)} 
  \to 0}
as $t\to\pm\infty$, 
where $H=\sqrt{-\Delta(2-\Delta)}$, and the wave operators $v_\pm\mapsto V^{-1}u(0)$ are local homeomorphisms around $0$ in $U^{-\si}H^{d/2-1}(\R^d)$. 
\end{theorem}
The corresponding result for the cubic NLS \eqref{NLS} holds \cite{CW} for any $d\ge 2$ and $|\si|\le d/2-1$. Moreover, the wave operators are defined on the whole energy space $H^1$ for $d=2,3,4$, and are bijective for $d=3,4$ \cite{GV,RyVi}. For $d=1$, modified wave operators have been constructed in \cite{O} for small data in weighted Sobolev spaces. 

The same result for NLS with the nonlinearity $F(u)$ is easily obtained for $d\ge 4$ and $|\si|\le d/2-2$ by the argument in \cite{CW}. 
The lower dimensional case is more difficult due to the quadratic terms and slow decay of the free evolution. For $d=3$, the wave operators have been constructed together with their inverses for the nonlinearities $|u|u$, $u^2$ and $\bar{u}^2$ for small data in weighted Sobolev spaces \cite{GOV,HMN}. For $d=2$, the modified wave operators for $|u|u$ and the ordinary ones for $u^2$ and $\bar{u}^2$ have been constructed for small data in weighted Sobolev spaces \cite{GO,HNST}. It has been proved that $|u|^2$ in $d=2$ does not admit the ordinary wave operators \cite{Shi}. 

The estimates in weighted spaces are quite sensitive to the linear part and so do not easily extend to our problem. Here the difficulty in our equation is that the Galilean invariance is also destroyed by the boundary condition. For the Maxwell-Higgs equation, a similar hyperbolic problem, the Lorentz invariance is preserved and the scattering result has been obtained in \cite{Ts}. 

To prove the above theorem, we first derive $L^p$ decay and Strichartz estimates for the linear evolution $e^{-iHt}$. These are slightly better for $d\ge 3$ around the Fourier origin than the free Schr\"odinger equation, because of curvature properties of the characteristic surfaces. This is the reason why we have a wider range of $\si$ than the NLS with the same nonlinearity. 
Actually, this gain is strong enough for $d\ge 5$ to cancel the $U^{-1}$ singularity in \eqref{eq v}, but not enough for $d\le 4$. Instead, we use a seemingly miraculous change of the unknown function -- a kind of ``low frequency normal form'' transformation -- which completely removes the singularities at the Fourier origin from the nonlinearity. 
However, it still leaves regular quadratic nonlinearities, whose decay is out of our control for $d\le 3$. 

In the next section 2, we investigate the linearized operator, deriving the decay and Strichartz estimates for the evolution. In Section 3, we perform the change of unknown function, and then in Section 4, prove the main theorem. 

\section{Linearized operator}
In this section, we derive the $L^p$ decay and Strichartz estimate for the linear equation 
\EQ{ \label{linear}
 i\dot u + \Delta u - 2\Re u = f.}
As in the introduction, we define the following $\C$-linear Fourier multipliers 
\EQ{
 U := \sqrt{-\Delta(2-\Delta)^{-1}},\quad H:=\sqrt{-\Delta(2-\Delta)},}
and $\R$-linear operator
\EQ{
 V u := U\Re u + i \Im u.}
By using the identities $H=-\Delta U^{-1}=(2-\Delta)U$, we obtain
\EQ{ \label{twist}
 i\dot u + \Delta u - 2 \Re u = iV(\p_t + iH)V^{-1}u.}
Thus the linear evolution by \eqref{linear} is reduced to the unitary group $e^{-iHt}$, for which we have the following decay and Strichartz estimates. Here $\dot B^s_p:=\dot B^s_{p,2}$ denotes the homogeneous Besov spaces (cf. \cite{BL}), and we denote spacetime norms by
$\| f \|_{L^{p} X} := \left( \int_{-\infty}^\infty \| f(t,\cdot) \|_X^p dt
\right)^{1/p}$, etc. 
\begin{theorem} \label{th:4-2}
Let $d\ge 1$. (i) For any $2\le q\le \infty$, we have
\EQ{
 \|e^{-itH}\fy\|_{\dot B^0_{q}} \lec t^{-d\si} \|U^{(d-2)\si}\fy\|_{\dot B^0_{q'}},
}
where $q'=q/(q-1)$ is the dual exponent and $\si=1/2-1/q$. 

(ii) For $j=1,2$, let $2\le p_j,q_j\le \infty$, $2/p_j + d/q_j = d/2$ and $s_j=\frac{d-2}{2}(1/2-1/q_j)$, but $(p_j,q_j)\not=(2,\infty)$. Then we have
\EQAL{
 &\|e^{-itH}\fy\|_{L^{p_1} \dot B^0_{q_1}} \lec \|U^{s_1}\fy\|_{L^2},\\
 &\left\|\int_{-\infty}^t e^{-i(t-s)H}f(s)ds\right\|_{L^{p_1} \dot B^0_{q_1}}
 \lec \|U^{s_1+s_2}f\|_{L^{p_2'} \dot B^0_{q_2'}}
}
\end{theorem}
When $d > 2$, these estimates are slightly better than those for the free Schr\"odinger equation because of the factor $U$, which is small around the Fourier origin. 
This gain is due to the fact that the curvature of the characteristic surface $|\xi|\sqrt{2+|\xi|^2}$ is larger than that of $|\xi|^2$ in the $(d-1)$ angular directions and smaller in the radial directions. 
The above theorem follows from the standard argument for the stationary phase. 
For completeness, we give a proof in the following subsection. 
\subsection{Stationary phase estimate}
\begin{theorem}
Let $\fy(r)\in C^\infty(0,\infty)$ satisfy the following. 
\begin{enumerate}
\item $\fy'(r),\fy''(r)>0$ for all $r>0$.
\item $\fy'(r)\sim \fy'(s)$ and $\fy''(r)\sim\fy''(s)$ for $0<s<r<2s$.
\item $|\fy^{(1+k)}(r)| \lec \fy'(r)/r^k$ for all $r>0$ and $k\in\N$.
\end{enumerate}
Let $\chi(r)$ be a dyadic cut-off function which is supported around $r\sim R$ and satisfies
\EQ{
 |\chi^{(k)}(r)| \lec R^{-k}.
}
(These estimates are supposed to hold uniformly for $r$ and $R$, but may depend on $k$.)
Then we have for any $d\in\N$,  
\EQ{
 \sup_{x\in\R^d} \left|\int_{\R^d} \chi(|\xi|) e^{i\fy(|\xi|)t + i\xi x} d\xi\right|
 \lec t^{-d/2} (\fy'(R)/R)^{-(d-1)/2} (\fy''(R))^{-1/2},
}
where the constant depends only on $\fy,\chi,d$.
\end{theorem}
\myproof
Since the case $d=1$ is easy, we may assume $d\ge 2$. The phase can be stationary when its gradient 
\EQ{
 \Phi := t \fy'(r)\th + x
}
vanishes, where we use the polar coordinates $\xi=r\th$. Because all other parts are radially symmetric, we may assume that $x=|x|e_1$, where $e_1=(1,0,\dots,0)\in\R^d$. 
By the monotonicity of $\fy'$, the stationary point for $\xi$ is uniquely determined for any fixed $(t,x)$, by $\xi=ce_1$ where $c$ is the unique solution (if it exists) for the equation
\EQ{
 |x| = t\fy'(c). 
}

First we dispose of the region away from the stationary point. 
There exists a smooth partition of unity on $S^{d-1}$, denoted by $F_1(\th),\dots F_N(\th) \in C^\infty(S^{d-1})$ where $N$ depends only on $d$, satisfying the following.
\begin{enumerate}
\item $\sum_j F_j(\th)=1$ for any $\th\in S^{d-1}$.
\item There exists $v^j\in S^{d-1}$ for each $j$ such that $\th\cdot v^j>1/2$ for any $\th\in\supp F_j$. 
\item For any $j>1$, $v^j \cdot x<0$.  
\end{enumerate}
Thus we have $\supp F_j\subset\{\th\cdot e_1\le\sqrt{3}/2\}$ for $j>1$. 
Denote by $I$ the integral to be estimated. 
We decompose the integral by using $F_j$ and estimate those terms for $j>1$ by partial integration in the $v^j$ direction. 
\EQ{
 I = \sum_{j=1}^N I_j, \quad
 I_j= \int \chi(r)F_j(\th) e^{-it\fy(r) + i\xi\cdot x} d\xi
}
\EQAL{
 I_j = \int e^{ip}(\p_j ip_j^{-1})^K(\chi(r)F_j(\th)) d\xi,
}
for any $K\ge 0$, where we denote
\EQ{
 p:= -t\fy(r) + \xi\cdot x, \quad \p_j=v_j\cdot\nabla, \quad p_j=\p_j p, 
}
and $\p_j ip_j^{-1}$ denotes the operator $f\mapsto \p_j (ip_j^{-1}f)$, not the derivative of $ip_j^{-1}$. This convention applies to similar expressions of the form $\p_* p_*^{-1}$ appearing later, but otherwise, we use $\p^\al$ to denote the usual partial derivatives. 

By the property of $F_j$ and $v^j$, we have
\EQ{
 p_j = v^j \cdot (-t\fy'(r)\th + x) < - t\fy'(r)/2,
}
for any $\th\in\supp F_j$. Expanding the $K$-th power, we estimate
\EQ{
 |I_j| \lec \sum_{|\al|=K} \int \left|\frac{\p_j^{\al_1}p_j \cdots \p_j^{\al_l} p_j}{p_j^{K+l}} \p_j^{\al_0}(\chi(r)F_j(\th)) \right| d\xi.
}
By using the estimates $|\fy^{(k+1)}|\lec|\fy'/r^k|$, $|\chi^{(k)}|\lec r^{-k}$ and $|\nabla^k F_j(\th)| \lec r^{-k}$, we obtain
\EQ{
 |I_j| \lec (t\fy'(R)R)^{-K}R^d
}
for any $j>1$ and $K\ge 0$. Choosing $K=d/2$, we obtain the desired estimate for those $I_j$, since $|\fy''|\lec|\fy'/r|$.  

Now we are left with the stationary part $I_1$. We further decompose this region into parallelepiped regions by a smooth partition of unity. Their side length should be $\de:=(t\fy''(R))^{-1/2}$ in the $e_1$ direction and $\la:=(t\fy'(R)/R)^{-1/2}$ in the other directions. 
Denote the set of such parallelepipeds by $\P$. 
For each $P\in\P$, we have a smooth cut-off function $G_P(\xi)$. We have $\sum_{P\in\P} G_P(\xi) = 1$ on $\supp F_1(\th)\chi(r)$, $\xi_1\sim |\xi| \sim R$ on any $P\in\P$, and 
\EQ{
 |D^\al G_P| \lec \de^{-\al_1}\la^{-\al_2-\dots-\al_n}.
}
Denote the integral corresponding to $G_P$ by $I_P$. 

We put the stationary point at the center of one parallelepiped $P_0$. In that region, the integral is not very oscillatory, so we just bound it by the volume:
\EQ{
 |I_{P_0}| \lec \de\la^{d-1},
}
which yields the desired estimate for this part. 

We have the same volume bound on the other parts, but we need additional decay factors to sum those numerous pieces up. For that purpose we may integrate by parts either in the $e_1$ direction or in the other $e^1\perp e_1$ directions. 

First we consider integration in $x^1$ direction. Take any $P\in\P$ on which we have $|\xi^1|\sim \la k$, $k\in\N$. Notice that such a $P$ exists only if $\la\lec R$. 
Then we have a direction $v\perp e_1$ such that $v\cdot\xi \sim \la k$ for any $\xi\in P$. We integrate by parts in the $v$ direction. Then we obtain
\EQ{
 I_P = \int e^{-it\fy(r) + i\xi\cdot x}(\p_v i p_v^{-1})^K (\chi(r)F_1(\th)G_P(\xi))d\xi ,
}
for any $K\ge 0$, where $\p_v=v\cdot\nabla$ and
\EQ{
 p_v=\p_v(-t\fy(r)+\xi\cdot x)= -t v\cdot\xi \fy'(r)/r \sim -t \la k \fy'(R)/R.
}
By expanding the $K$-th power, we obtain
\EQ{
 |I_P| \lec (t\la k \fy'(R)R \min(\la,R))^{-K} \la^{d-1}\de,
}
where $\la$ in $\min(\la,R)$ is coming from derivatives hitting $G_P(\xi)$. 
Thus we get 
\EQ{ \label{bound in e^1}
 |I_P| \lec k^{-K} \de\la^{d-1}. 
}

Next we consider integration in the $e_1$ direction. 
Here we should be a bit more careful, since the phase can be very small depending on $\fy''$ even if we are not in $P_0$, and then the estimate on the higher derivatives of the phase becomes less trivial. 
Take any $P\in\P$ on which we have $\xi_1-c\sim \de j$ with $j\not=0$, where $c$ denotes the $\xi_1$ value of the stationary point. Such a $P$ exists only if $\de\lec R$. 
Let $p_1=\p_1 p=-t\fy'(r)\th_1+|x|$. We compute 
\EQAL{
 &\p_1 r = \th_1, \quad \p_1\th_1 = (1-\th_1^2)/r,\\ 
 &\p_1 p_1 = -t(\fy''\th_1 + \fy'(1-\th_1^2)/r)
}
For higher derivatives we have 
\EQ{
 |\p_1^{k+1} p_1| \lec \sum_{j=0}^k t |\fy^{(2+j)}/r^j| + t|\fy'(1-\th_1^2)/r^{k+1}|.
}
Since $\fy',\fy''>0$, we have $|\p_1^{k+1} p_1|\lec |\p_1 p_1|/r^k$. This in particular implies that $\p_1 p_1$ is of the same size on $\supp\chi$. 
Then it implies that $|p_1|\sim|\p_1 p_1 \de j|$. We have also $|\p_1 p_1|\gec t|\fy''(R)|$. 

Integrating by parts in the $e_1$ direction we have
\EQ{
 I_P = \int e^{-it\fy(r) + i\xi\cdot x}(\p_1 i p_1^{-1})^K (\chi(r)F_1(\th)G_P(\xi))d\xi \
,
}
Expanding the $K$-th power, we estimate
\EQ{
 |I_P| \lec (t\fy''(R)\de j \de)^{-K} \la^{d-1}\de,
}
Thus we obtain 
\EQ{
 |I_P| \lec j^{-K} \la^{d-1}\de. 
}
In conclusion, we have for any $P\in \P$,
\EQ{
 |I_P| \lec (1+k+j)^{-K} \la^{d-1}\de,
}
for any $K\in\N$, so we can sum them up for all $k,j$, deriving the desired estimate.  
\myendproof

\section{Low frequency normal form}
Let $u$ be a solution of \eqref{GLS2} and $v=V^{-1}u$. 
By using \eqref{twist}, we can derive the equation for $v$: 
\EQ{
  v(t) = e^{-itH}v_0 - \int_0^t e^{-i(t-s)H} V^{-1}iF(u(s))\, ds.}
The presence of $U^{-1}$ in the nonlinearity is our only ``enemy'' for $d\ge 4$; without it, the scattering would follow in the same way as for NLS with quadratic and cubic power nonlinearities. 
Quadratic nonlinearities have the critical decay rate for scattering by the Strichartz estimate when $d=4$, while for $d>4$ their decay is more than sufficient. Indeed we can derive the scattering result for $d>4$ and certain restricted values of $\si$ by using only the Strichartz estimate with the $U$ gain and the Sobolev embedding. That argument, however, cannot work in the critical case $d=4$. 

We introduce the following nonlinear change of variable to kill the singular nonlinearity. It gives a better range of $\si$ even for $d>4$. Let 
\EQ{ \label{def w}
 w := N(u) := u + P|u|^2/2,}
where $P$ is the Fourier multiplier $F^{-1}\chi\F$, $\chi\in C_0^\infty(\R^d)$ satisfying $\chi(x)=1$ for $|x|\le 1$ and $\chi(x)=0$ for $|x|\ge 2$. 
Then $w$ satisfies the following equation:
\EQAL{
 &i\dot w= -\Delta w + 2\Re w + G_1 + i G_2,\\
 &G_1 = (3-P)u_1^2 + Q u_2^2 + P\Delta|u|^2/2 + |u|^2 u_1,\\ 
 &G_2 = 2Q(u_1 u_2) + \nabla P\cdot(u_2\nabla u_1 - u_1\nabla u_2) +Q(|u|^2u_2),}
where $u=u_1+iu_2$ and $Q=Id-P$. Thus 
\[
  z:=V^{-1}w=V^{-1}N(u)=U^{-1}(u_1+P|u|^2/2)+iu_2=v+U^{-1}P|Vv|^2/2=:M v
\]
satisfies
\EQ{
 z=e^{-iHt}z(0) - \int_0^t e^{-iH(t-s)}(iG_1-U^{-1}G_2) ds.}
Notice that the singularity of $U^{-1}$ in front of $G_2$ is now canceled by either $Q$ or $\nabla$. There is no loss of regularity thanks to $P$. Moreover, $G$ does not contain bilinear interaction of very low frequency of $u_2$, whereas $u_1$ behaves better at the Fourier $0$ due to the relation with $z$. 

\section{Proof of the main theorem}
We introduce Besov spaces with different regularity for low and high frequency:
\EQAL{
 \|\fy\|_{B^{a,b}_q} := \|P\fy\|_{\dot B^a_q} + \|Q\fy\|_{\dot B^b_q},\quad 
 H^{a,b}:=B^{a,b}_2,}
where $\dot B^a_q:=\dot B^a_{q,2}$ denotes the homogeneous Besov space (cf. \cite{BL}). We will use these spaces under the restriction $q<d/a$, embedding them into $\S'$. We have $B^{a,b}_q = U^{b-a}\dot B^b_{q,2}$, $H^{0,b}=H^b$ and the embedding $B^{a,b}_q \subset B^{a',b'}_q$ for $a\le a'$ and $b\ge b'$. 
These spaces have the following multiplicative property, which is enough for us to estimate the nonlinearity.
\begin{lemma}
Let $p_1,p_2,p_3\in[2,\infty]$, $s_1,s_2,s_3\in\R$, $t_1,t_2,t_3\in\R$ and $b_a:=1/p_a$ for $a=1,2,3$. Assume that for $a=1,2,3$, 
\EQAL{ \label{cond exp}
 \max(0,s_a,s_1+s_2+s_3) \le d(b_1+b_2+b_3-1) \le t_1+t_2+t_3\ge t_a,\ s_a < d b_a.}
Then we have
\EQ{
 \left|\int_{\R^d} f g h dx\right| \lec \|f\|_{B^{s_1,t_1}_{p_1}}\|g\|_{B^{s_2,t_2}_{p_2}}\|h\|_{B^{s_3,t_3}_{p_3}}.}
\end{lemma}
\myproof
We apply the Littlewood-Paley decomposition to each of the functions. 
\EQ{
 \int fgh = \sum_{l=-\infty}^\infty \sum_{k = l-2}^l \sum_{j=-\infty}^k c_{j,k,l} \int f_j g_k h_l + f_k g_l h_j + f_l g_j h_k,}
where $0\le c_{j,k,l}\le 1$ and each $f_j$ is supported in Fourier space in $\{2^{j-1}\le|\xi|\le 2^{j+1}\}$. By symmetry, it suffices to estimate the first summand. Let $\be=1-b_2-b_3$, then we have $0\le\be\le b_1$. By the Sobolev and the H\"older inequalities and the embedding $\dot B^0_{p}\subset L^{p}$, we have
\EQAL{
 \biggl|\sum_{j\le k} \int f_j g_k h_l\biggr| &\le \|\sum_{j\le k}f_j\|_{L^{1/\be}} \|g_k\|_{L^{p_2}} \|h_l\|_{L^{p_3}}\\
 &\lec 
   \|\sum_{j\le k} 2^{j(-s+d[b_1-\be])} f_j\|_{\dot B^s_{p_1}} \|g_k\|_{L^{p_2}} \|h_l\|_{L^{p_3}}}
for any $s$. Taking $s=s_1$ for $j \leq 0$ and $s=t_1$ for $j>0$ (when $k>0$),
the norm on $f$ is bounded by 
\EQ{ \label{sum in j}
 \begin{cases}
 2^{k(-s_1+d(b_1-\be))} &(k\le 0),\\
 1 + 2^{k(-t_1+d(b_1-\be))} &(k> 0), 
 \end{cases}}
since $-s_1+d(b_1-\be)\ge 0$. 
The summation over $k \sim l\in\Z$ converges, because for $k\le 0$,  
\EQAL{
 -s_1+d(b_1-\be) - s_2 - s_3
 = -(s_1+s_2+s_3) + d(b_1+b_2+b_3 - 1) \ge 0.}
and for $k\ge 0$, 
\EQ{
 -t_2-t_3 \le 0,\quad -t_1-t_2-t_3 + d(b_1+b_2+b_3-1) \le 0,}
and the summability is provided by the $\ell^2$ norm in the Besov norms.
\myendproof

For the proof of our main theorem, we will use the following special cases of the above lemma. Let $d\ge 4$, $s=d/2-1$, $b=1/p=1-1/p'=1/2-1/d$ and $1/q=1/2-1/(2d)$. Then we have 
\EQAL{ \label{bil}
 &(1)\quad H^{\si-b,s-j} \times B^{\si-b,s-k}_p \subset B^{\si+b,s-j-k}_{p'},\\
 &(2)\quad H^{\si,s} \times B^{\si-b,s-1}_p \subset H^{\si-b,s-1},\\ 
 &(3)\quad B^{\si-b/2,s}_q \times B^{\si-b/2,s}_q \subset H^{\si,s},\\
 &(4)\quad H^{\si,s} \times H^{\si,s} \subset B^{\si,s}_{p'},\\
 &(5)\quad H^{\si,s} \times B^{\si-b/2,s}_q \subset H^{\si-1,s-1/2},
}
for $0\le j,k,j+k\le 1$, provided that $|\si|-b\le d/2-2$. All conditions in the lemma are satisfied because
\EQAL{
 (1)\quad &\max(\si-3b,|\si|-b)\le d/2-2\le s\ge j+k-s,\ |\si|-b<d/2-1,\\[4pt]
 (2)\quad &\max(\si,b-\si) \le d/2-1\le s\ge 1-s,\ \si<d/2,\\[4pt]
 (3)\quad &\max(\si-b/2, -\si) \le d/2-1\le s,\ \si-b/2<d/2-1/2,\\[4pt]
 (4)\quad &|\si|\le d/2-1\le s,\\[4pt]
 (5)\quad &\max(\si-b/2+1,\si,1-\si)\le d/2-1/2\le s+1/2\ge 1/2-s.}
As the function space for $v=V^{-1}u$ and $z=V^{-1}N(u)$, we define $X_0$ and $X_1$ by
\EQ{
 X_j = L^\infty H^{\si-jb,s} \cap L^2 B^{\si-b,s}_p.}
We have $U:X_0\to X_1\subset X_0$ bounded. First we consider the norm relation between $v$ and $z$. By \eqref{bil}-(2)(4) and the embeddings $B^{\si,s}_{p'}\subset H^{\si-1,s-1}$, $H^{\si,s}\subset B^{\si-b,s-1}_p$ and $H^{\si-b,s-1}\subset B^{\si-b,s-2}_p$, we have
\EQ{ \label{bil3}
 \|U^{-1}P(fg)\|_{B} \lec \|f\|_{H^{\si,s}} \|g\|_{B}}
for $B=H^{\si,s}$ and $B=B^{\si-b,s}_p$. 
Hence for any small $z\in H^{\si,s}$, the unique inverse $v=M^{-1}z$ is obtained in the same space by applying the Banach fixed point theorem to the equation $v=z-U^{-1}P|Vv|^2/2$. Thus $M$ is bi-Lipschitzian for small data in $H^{\si,s}$. In particular, we have the following estimate for the free part $z_0:=e^{-iHt}z(0)$
\EQ{ \label{free est}
 \|z_0\|_{X_0} \lec \|z(0)\|_{H^{\si,s}} \lec \|v(0)\|_{H^{\si,s}}<\de.}
Moreover, \eqref{bil3} implies that $\|z\|_{X_0}\sim \|v\|_{X_0}$ as long as they remain small. In addition, we have
\EQ{ \label{M triv}
 \|M(e^{-iHt}\fy)-e^{-iHt}\fy\|_{H^{\si,s}} \to 0 \quad(|t|\to\infty)}
for any small $\fy\in H^{\si,s}$. This is proved as follows. First if $\fy\in C_0^\infty$, then the decay estimate and \eqref{bil}-(5) implies that
\EQ{
 \|U^{-1}P|Ve^{-iHt}\fy|^2\|_{H^{\si,s}} \lec \|e^{-iHt}\fy\|_{H^{\si,s}} \|e^{-iHt}\fy\|_{B^{\si-b/2,s}_q} \lec t^{-1/2}.}
Then \eqref{M triv} for general $\fy$ follows from the denseness of $C_0^\infty\subset H^{\si,s}$ and the Lipschitz continuity of $M$. The bi-Lipschitz property implies also that
\EQ{ \label{M-1 triv}
 \|e^{-iHt}\fy-M^{-1}(e^{-iHt}\fy)\|_{H^{\si,s}} \to 0 \quad(|t|\to\infty)}
for any small $\fy\in H^{\si,s}$. 

Next we estimate the quadratic terms in $G(u)$. By \eqref{bil}-(1), we have
\EQAL{
 &\|u_1u\|_{L^2 B^{\si+b,s}_{p'}} \lec \|u_1\|_{L^\infty H^{\si-b,s}} \|u\|_{L^2 B^{\si-b,s}_p} \lec \|u_1\|_{X_1} \|u\|_{X_0},\\
 &\|P(u\na u)\|_{L^2 B^{\si+b,s}_{p'}} \lec \|u\|_{L^2 B^{\si-b,s}_p} \|\na u\|_{L^\infty H^{\si-b,s-1}} \lec \|u\|_{X_0}^2,}
These are used to estimate the terms $u_1^2$, $u_1u_2$, $P\Delta|u|^2$ and $u_2\na u_1-u_1\na u_2$, together with the boundedness of $U^{-1}Q$, $U^{-1}\na P$ and $P\na$ on every Besov space. $Qu_2^2$ can be treated in the same way as $u_1 u$ since 
\EQ{
 Q u_2^2 = Q(u_2 \cdot Q_{-2} u_2) + Q(Q_{-2} u_2 \cdot P_{-2} u_2)}
and $Q_{-2}:X_0\to X_1$ bounded, where $P_{-2}:=\F^{-1}\chi(4\xi)\F$ and $Q_{-2}:=Id-P_{-2}$, defined in the same way as $P$ and $Q$ (see below \eqref{def w}).  

For the cubic terms, we use \eqref{bil}-(3)(4) together with the interpolation inequality $X_0\subset L^4 B^{\si-b/2,s}_q$, yielding
\EQAL{
 \|u^3\|_{L^2 B^{\si,s}_{p'}}
 \lec \|u\|_{L^\infty H^{\si,s}} \|u^2\|_{L^2 H^{\si,s}}
 \lec \|u\|_{X_0} \|u\|_{L^4 B^{\si-b/2,s}_q}^2 \le \|u\|_{X_0}^3.}

Gathering these estimates, we obtain
\EQAL{
 \|z-z_0\|_{X_0} &\lec \|V^{-1}iG\|_{L^2 B^{\si+b,s}_{p'}} \lec \|u\|_{X_0}(\|u\|_{X_0}+\|u_1\|_{X_1}+\|u\|_{X_0}^2) \lec \|z\|_{X_0}^2,}
if the right hand side is small. By continuity in time and \eqref{free est}, we conclude that
\EQ{
 \|z\|_{X_0} + \|z_0\|_{X_0} + \|V^{-1}iG\|_{L^2 B^{\si+b,s}_{p'}} \lec \|v(0)\|_{H^{\si,s}}.}
In particular, the norm of $G$ on the interval $(T,\infty)$ vanishes as $T\to\infty$, which implies via the Strichartz estimate that $e^{iHt}z(t)\to\exists v_\pm$ strongly in $H^{\si,s}$. Then \eqref{M-1 triv} implies that $e^{iHt}v(t)\to v_\pm$ in $H^{\si,s}$. 

We can construct the wave operators applying the same estimates to the Cauchy problem from $t=\pm\infty$:
\EQ{
 z = e^{-iHt}v_\pm - \int_{\pm\infty}^t e^{-iH(t-s)}V^{-1}iG(VM^{-1} z(s)) ds.}
The bi-Lipschitz property of the maps $v(0)\mapsto v_\pm$ is similarly proved by estimating the difference of two solutions in the same space. These arguments, essentially by the Banach fixed point theorem, are quite standard in nonlinear scattering theory, and so we omit the details. 

\section*{Acknowledgments}
The research of Gustafson and Tsai is partly supported by NSERC grants
nos.~22R80976 and 22R81253. The research of Nakanishi is partly
supported by the JSPS grant no.~15740086.

\noindent{Stephen Gustafson},  gustaf@math.ubc.ca \\
Dept. Mathematics, University of British Columbia, 
Vancouver, BC V6T 1Z2, Canada

\bigskip

\noindent{Kenji Nakanishi},
n-kenji@math.kyoto-u.ac.jp \\
Dept. Mathematics, Kyoto University,
Kyoto 606-8502, Japan

\bigskip

\noindent{Tai-Peng Tsai},  ttsai@math.ubc.ca \\
Dept. Mathematics, University of British Columbia, 
Vancouver, BC V6T 1Z2, Canada


\begin{thebibliography}{10}


\bibitem{BL}
J.~Bergh and J.~{L\"{o}fstr\"{o}m}, 
{Interpolation spaces}, An introduction. 
Grundlehren Math. Wiss. {\bf 223}, Springer,
Berlin--Heiderberg--New York, 1976.

\bibitem{BOS} F. Bethuel, G. Orlandi, D. Smets,
{\it Vortex rings for the Gross-Pitaevskii equation}.
J. Euro. Math. Soc. {\bf 6} (2004), no. 1, 17--94.

\bibitem{BS} F.~Bethuel and J.~C.~Saut, 
{\it Travelling waves for the Gross-Pitaevskii equation. I}. 
Ann. Inst. H. Poincar\'e Phys. Th\'eor. {\bf 70} (1999), no. 2, 147--238.
 
\bibitem{Ch} D.~Chiron, 
{\it Travelling waves for the Gross-Pitaevskii equation in dimension
  larger than two}. 
Nonlinear Anal. {\bf 58} (2004), no. 1-2, 175--204. 

\bibitem{CW} T.~Cazenave and F.~B.~Weissler, 
{\it The Cauchy problem for the critical nonlinear Schr\"odinger
  equation in $H\sp s$}. 
Nonlinear Anal. {\bf 14} (1990), no. 10, 807--836. 

\bibitem{CJ} J.~E.~Colliander and R.~L.~Jerrard, 
{\it Vortex dynamics for the Ginzburg-Landau-Schr\"odinger equation}. 
Internat. Math. Res. Notices {\bf 1998}, no. 7, 333--358; and
{\it Ginzburg-Landau vortices: weak stability and Schr\"odinger
  equation dynamics}. 
J. Anal. Math. {\bf 77} (1999), 129--205. 

\bibitem{FS} A.L. Fetter and A.A. Svidzinsky,
{\it Vortices in a trapped dilute Bose-Einstein condensate}.
Preprint: arXiv:cond-mat/0102003.

\bibitem{GO} J.~Ginibre and T.~Ozawa, 
{\it Long range scattering for nonlinear Schr\"odinger and Hartree
  equations in space dimension $n\geq 2$}. 
Comm. Math. Phys. {\bf 151} (1993), no. 3, 619--645.
 
\bibitem{GOV} J.~Ginibre, T.~Ozawa and G.~Velo, 
{\it On the existence of the wave operators for a class of nonlinear
  Schr\"odinger equations}. 
Ann. Inst. H. Poincar\'e Phys. Th\'eor. {\bf 60} (1994), no. 2, 211--239. 

\bibitem{GV} J.~Ginibre and G.~Velo, 
{\it Scattering theory in the energy space for a class of nonlinear
  Schr\"odinger equations}. 
J. Math. Pures Appl. (9) {\bf 64} (1985), no. 4, 363--401.

\bibitem{Gr} P.~Gravejat, 
{\it Limit at infinity for travelling waves in the Gross-Pitaevskii equation}. 
C. R. Math. Acad. Sci. Paris {\bf 336} (2003), no. 2, 147--152; and 
{\it A non-existence result for supersonic travelling waves in the
  Gross-Pitaevskii equation}. 
Comm. Math. Phys. {\bf 243} (2003), no. 1, 93--103. 

\bibitem{Gu} S.~Gustafson, 
{\it Dynamic stability of magnetic vortices}. 
Nonlinearity {\bf 15} (2002), no. 5, 1717--1728. 

\bibitem{GS} S. Gustafson, I.M. Sigal,
{\it Effective dynamics of magnetic vortices.}
To appear in Adv. Math. (2005).

\bibitem{HMN} N.~Hayashi, T.~Mizumachi and P.~I.~Naumkin, 
{\it Time decay of small solutions to quadratic nonlinear
  Schr\"odinger equations in 3D}. 
Differential Integral Equations {\bf 16} (2003), no. 2, 159--179.

\bibitem{HNST} N.~Hayashi, P.~I.~Naumkin, A.~Shimomura and
  S.~Tonegawa, 
{\it Modified wave operators for nonlinear Schr\"odinger equations in
  one and two dimensions}. 
Electron. J. Differential Equations {\bf 2004}, No. 62, 16 pp.

\bibitem{J} R. Jerrard,
{\it Vortex dynamics for the Ginzburg-Landau wave equation}.
Calc. Var. Partial Diff. Eqns. {\bf 9} (1999) no.8, 683-688. 

\bibitem{J2} R.L. Jerrard,
{\it Vortex filament dynamics for Gross-Pitaevsky type equations.}
Ann. Sc. Norm. Super. Pisa Cl. Sci. {\bf 1} (2002) no. 4, 733-768.


\bibitem{JPR} C.A. Jones, S.J. Putterman, P.H. Roberts,
{\it Motions in a Bose condensate: stability of solitary
wave solutions of non-linear Schr\"odinger equations
in two and three dimensions}.
J. Phys. A: Math. Gen. {\bf 19} (1986) 2991-3011.

\bibitem{LS} O.~Lange and B.~J.~Schroers, 
{\it Unstable manifolds and Schr\"odinger dynamics of Ginzburg-Landau vortices}. Nonlinearity {\bf 15} (2002), no. 5, 1471--1488.

\bibitem{L} F.-H. Lin,
{\it Vortex dynamics for the nonlinear wave equation.}
Comm. Pure Appl. Math. {\bf 52} (1999) no.6, 737-429.

\bibitem{LX} F.~H.~Lin and J.~X.~Xin, 
{\it On the incompressible fluid limit and the vortex motion law of
  the nonlinear Schr\"odinger equation}. 
Comm. Math. Phys. {\bf 200} (1999), no. 2, 249--274. 

\bibitem{OS} Y.~N.~Ovchinnikov and I.~M.~Sigal, 
{\it Long-time behaviour of Ginzburg-Landau vortices}. 
Nonlinearity {\bf 11} (1998), no. 5, 1295--1309. 

\bibitem{O} T.~Ozawa, 
{\it Long range scattering for nonlinear Schr\"odinger equations in
  one space dimension}. 
Comm. Math. Phys. {\bf 139} (1991), no. 3, 479--493.

\bibitem{RyVi} E.~Ryckman and M.~Visan, 
{\it Global well-posedness and scattering for the defocusing
  energy-critical nonlinear Schr\"odinger equation in $\R^{1+4}$}, 
preprint, http://jp.arxiv.org/abs/math.AP/0501462.

\bibitem{Shi} A.~Shimomura, 
{\it Nonexistence of asymptotically free solutions for quadratic
  nonlinear Schr\"odinger equations in two space dimensions}. 
Differential Integral Equations {\bf 18} (2005), no. 3, 325--335. 

\bibitem{Sp} D.~Spirn, 
{\it Vortex motion law for the Schr\"odinger-Ginzburg-Landau
  equations}. 
SIAM J. Math. Anal. {\bf 34} (2003), no. 6, 1435--1476.

\bibitem{Ts} Y.~Tsutsumi, 
{\it Stability of constant equilibrium for the Maxwell-Higgs
  equations}. 
Funkcial. Ekvac. {\bf 46} (2003), no. 1, 41--62. 

\end{thebibliography}
\end{document}